\documentclass[reqno]{amsart}
\usepackage{amsmath}

\def\eqn#1{(\ref{eq:#1})}
\newcommand{\qBin}[3]{\genfrac{[}{]}{0pt}{0}{#1}{#2}_{#3}}
\newcommand{\qBinsm}[3]{\genfrac{[}{]}{0pt}{1}{#1}{#2}_{#3}}
\DeclareMathOperator{\mymod}{mod}
\theoremstyle{plain}
\newtheorem*{lem}{Lemma}

\numberwithin{equation}{section}
\allowdisplaybreaks[3]

\begin{document}

\title{Variants of the Andrews-Gordon Identities}
\author[A.~Berkovich]{Alexander Berkovich}
\address{Department of Mathematics, The Pennsylvania State University,
         University Park, PA~16802, USA}
\email{alexb@math.psu.edu}
\thanks{The first author was partially supported by SFB-grant
        F1305 of the Austrian FWF}
\author[P.~Paule]{Peter Paule}
\address{Research Institute for Symbolic Computation,
         Johannes Kepler University,
	 A--4040 Linz, Austria}
\email{Peter.Paule@risc.uni-linz.ac.at}

\begin{abstract}
The object of this paper is to propose and prove a new generalization 
of the Andrews-Gordon Identities, extending a recent result of Garrett, Ismail and 
Stanton. We also give a combinatorial discussion of the finite form of their result
which appeared in the work of Andrews, Knopfmacher, and Paule.
\end{abstract}

\maketitle

\section{Introduction}\label{sec:1}

The celebrated Rogers-Ramanujan identities are given analytically as follows
\begin{align}
\sum_{t\ge0} \frac{q^{t^2+at}} {(q)_t} & =
\frac{1}{(q)_\infty} \sum_{j=-\infty}^\infty
\left\{q^{j(10j+1+2a)}-q^{(2j+1)(5j+2-a)}\right\}\notag\\
 & = \frac{1} {(q^{1+a},q^5)_\infty(q^{4-a},q^5)_\infty},
\label{eq:1.1}
\end{align}
where $a=0,1$ and the $q$-shifted factorials $(z;q)_t$ are defined as 
usual as
\begin{equation}
(z;q)_t=(z)_t=\begin{cases}
\prod_{j=0}^{t-1}(1-zq^j), & \text{if } t\in \mathbb{Z}_{>0}, \\
1, & \text{if } t=0. \end{cases}
\label{eq:1.2}
\end{equation}
It is well known that these identities have polynomial analogs. 
In particular, building on the work of Schur and MacMahon,
Andrews~\cite{A1} has shown that for $L\in \mathbb{Z}_{\geq 0}$ 
\begin{equation}
\sum_{t\geq0} q^{t^2}\qBin{L-t}{t}{q}= e_L(q) 
\label{eq:1.3}
\end{equation}
and
\begin{equation}
\sum_{t\ge0}q^{t^2+t} \qBin{L-t-1}{t}{q} = d_L(q), 
\label{eq:1.4}
\end{equation}
where
\begin{equation}
e_L(q)=\sum_{j=-\infty}^\infty \left\{ q^{j(10j+1)}
\qBin{L}{\lfloor\frac{L}{2}\rfloor-5j}{q}
-q^{(2j+1)(5j+2)} \qBin{L}{\lfloor\frac{L-4}{2}\rfloor-5j}{q}
\right\}
\label{eq:1.5}
\end{equation}
and
\begin{equation}
d_L(q)=\sum_{j=-\infty}^\infty \left\{ q^{j(10j+3)}
\qBin{L}{\lfloor\frac{L-1}{2}\rfloor-5j}{q}
-q^{(2j+1)(5j+1)}
\qBin{L}{\lfloor\frac{L-3}{2}\rfloor-5j}{q}
\right\}.
\label{eq:1.6}
\end{equation}
As usual, $\lfloor x \rfloor$ denotes the integer part of $x$ and $q$-binomial
coefficients are defined as follows
\begin{equation}
\qBin{n+m}{n}{q} =
\begin{cases}\frac{(q^{n+1})_m}{(q)_m}, & \mbox{if } m\in\mathbb{Z}_{\ge0}, \\
0, & \mbox{otherwise.} \end{cases}
\label{eq:1.7}
\end{equation}
Both polynomial sequences $(e_L)$ and $(d_L)$ satisfy the recurrence
\begin{equation}
c_L(q)=c_{L-1}(q)+q^{L-1}c_{L-2}(q), \quad L\ge 2.
\label{eq:1.8}
\end{equation}
The above equation along with the initial conditions
\begin{equation}
d_0(q)=0,\quad e_0(q)=e_1(q)=d_1(q)=1
\label{eq:1.9}
\end{equation}
specifies these sequences uniquely. Moreover, one can read \eqn{1.8} backward to 
define $(e_L)$, $(d_L)$ for negative subindices; i.e., for $L \geq 1$,
\begin{equation}
\begin{cases}
e_{-L}(q)=(-1)^L q^{\binom{L}{2}}d_{L-1}(\frac{1}{q}), \\
d_{-L}(q)=(-1)^{L+1} q^{\binom{L}{2}}e_{L-1}(\frac{1}{q}).\end{cases}
\label{eq:1.10}
\end{equation}
Despite the long history of the Rogers-Ramanujan identities, the 
following variants found by Garrett et al.~\cite{GIS}
\begin{equation}
\sum_{t\ge0} \frac{q^{t^2+mt}} {(q)_t}=
\frac{(-1)^m q^{-\binom{m}{2}}d_{m-1}(q)}{(q,q^5)_\infty(q^4,q^5)_\infty}+                    
\frac{(-1)^{m+1}q^{-\binom{m}{2}}e_{m-1}(q)}
{(q^2,q^5)_\infty(q^3,q^5)_\infty},\;\;\;\;\; m\ge 0
\label{eq:1.11}
\end{equation}
appeared to be new, even though closely related results were derived before 
in \cite{C} and \cite{AB}.

Actually, \eqn{1.11} can be extended to negative $m$ with the aid of \eqn{1.10} as
\begin{equation}
\sum_{t\ge 0} \frac{q^{t^2-Mt}} {(q)_t}=
\frac{e_M (\frac{1}{q})}{(q,q^5)_\infty(q^4,q^5)_\infty}+
\frac{d_M (\frac{1}{q})}{(q^2,q^5)_\infty(q^3,q^5)_\infty}
\label{eq:1.12}
\end{equation}
with $M\in\mathbb{Z}_{\ge0}$. The authors of \cite{GIS} gave two proofs of \eqn{1.11}.
In the first proof they evaluated a certain integral involving $q$-Hermite polynomials
in two different ways and equated the results. Their second proof made essential 
use of Schur's involution. A very different approach was taken by Andrews et al.\
in~\cite{AKP},
where identity \eqn{1.12} appeared as a limiting case of the much stronger identity 
\begin{align}
\sum_{t\ge0} q^{t^2+mt}\qBin{L-t}{t}{q} &=
(-1)^m q^{-\binom{m}{2}}d_{m-1}(q)e_{L+m}(q)\notag\\
&+ (-1)^{m+1}q^{-\binom{m}{2}}e_{m-1}(q)d_{L+m}(q), \quad
(L,m\ge0),
\label{eq:1.13}
\end{align}
which was proven recursively. It is trivial to verify that in the limit $L\rightarrow\infty$ 
\eqn{1.13}
turns into \eqn{1.11}. It was pointed out in \cite{AKP} that \eqn{1.13} may be viewed as a
$q$-analog of a famous Euler-Cassini's identity for Fibonacci
numbers. We remark that a new approach to identities of
$q$-Euler-Cassini type has been given in \cite{Cigler}.

Once again, we
may employ \eqn{1.10} to extend \eqn{1.13} to negative $m$. The result is
\begin{equation}
\sum_{t\ge 0}q^{t^2-Mt}\qBin{L-t}{t}{q}=
e_M(\frac{1}{q})e_{L-M}(q)+d_M(\frac{1}{q})d_{L-M}(q), \quad M\ge 0.
\label{eq:1.14}
\end{equation}
Remarkably, this reformulation of \eqn{1.13} enables us to reduce it to \eqn{1.3}, 
\eqn{1.4} in an elementary combinatorial fashion. This is done in Section~\ref{sec:2}.
In Section~\ref{sec:3}, we briefly discuss a polynomial version of the Andrews-Gordon identities
and then, move on to our main results \eqn{3.19}--\eqn{3.21}: variants of the Andrews-Gordon
identities, which are straightforward multisum generalizations of \eqn{1.14}.
In Section~\ref{sec:4}, some problems for further investigation motivated by this work are indicated.
Finally, certain technical details pertaining to the recurrences for multisums are relegated 
to the Appendix.

\section{Combinatorial analysis of \eqn{1.14}}\label{sec:2}

We start by recalling a well-known fact.

\begin{lem}
For $L\ge 2t\ge0$, 
$q^{t^2}\qBinsm{L-t}{t}{q}$ 
is the generating function for partitions into exactly $t$ parts with 
difference at least $2$ between parts, such that each part $< L$.
\end{lem}

We wish to describe this generating function in ``path" language. To this end we define an 
admissible sequence of integers $\Sigma$ as an ordered sequence 
$(\sigma_i,\sigma_{i+1},\dots,\linebreak[1]\sigma_{f-1},\sigma_f)$ such that $\sigma_l \in \{0,1\}$ 
for $i\le l\le f$ and $\sigma_j\sigma_{j+1}=0$ for $i\le j\le f-1$. Given $\Sigma$ we can 
construct an admissible lattice path $P(\Sigma)$ by connecting points $(j;\sigma_j)$ and 
$(j+1;\sigma_{j+1})$ by the straight line segments. Thus, any admissible path is made out 
of three basic segments:
\begin{center}
\begin{picture}(80,15)(200,0)
\linethickness{0.7pt}
\put(200,0){\line(3,2){15}}
\put(216,0){,}
\put(240,0){\line(-3,2){15}}
\put(245,0){,}
\put(255,0){\line(2,0){20}}
\put(280,0){.}
\end{picture}
\end{center}

\begin{figure}[ht]
\begin{center}
\setlength{\unitlength}{1cm}
\begin{picture}(12.5,3.2)(-.25,-.4)
\linethickness{0.2pt}
\put(0,0){\line(0,1){2}}
\put(0,0){\line(1,0){12}}
\put(0,2){\vector(0,1){0}}
\put(12,0){\vector(1,0){0}}
\put(0,1){\circle*{.10}}
\multiput(0,0)(1,0){11}{\circle*{.10}}
\put(-.25,-.4){\mbox{$0$}}
\put(1.9,-.4){\mbox{$2$}}
\put(3.9,-.4){\mbox{$4$}}
\put(5.9,-.4){\mbox{$6$}}
\put(7.9,-.4){\mbox{$8$}}
\put(9.9,-.4){\mbox{$L=10$}}
\put(12.2,-.4){\mbox{$j$}}
\put(-.25,.9){\mbox{$1$}}
\put(-.25,2.2){\mbox{$\sigma_j$}}
\thicklines
\put(0,0){\line(1,0){1}}
\put(1,0){\line(1,1){1}}
\put(2,1){\line(1,-1){1}}
\put(3,0){\line(1,1){1}}
\put(4,1){\line(1,-1){1}}
\put(5,0){\line(1,0){2}}
\put(7,0){\line(1,1){1}}
\put(8,1){\line(1,-1){1}}
\put(9,0){\line(1,0){1}}
\end{picture}
\parbox{0.8\textwidth}{\caption{\label{fig:AXEL}Path representation of
the partition $14=2+4+8$ into $3$ parts, each $< L = 10$.}}
\end{center}
\end{figure}
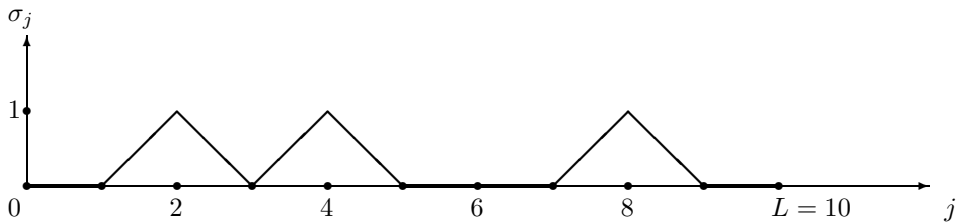

\noindent
Note that a horizontal segment is always of height $0$.
On such a path we distinguish points $(j;\sigma_j=1)$ 
with $i \neq j$ and $i\neq f$,
which we call peaks. Clearly, 
the distance between two peaks is at least $2$, as can be seen from Figure~\ref{fig:AXEL}.

Let us denote the space of all admissible paths $P(\sigma_i=s,\sigma_{i+1},\dots,\sigma_f=b)$ 
with exactly $t$ peaks and fixed end points $(i;s),(f;b)$ as $\mathbf{P}_{s,b}^t(i,f)$. 
For a given path $\in \mathbf{P}_{0,0}^t(0,L)$ we can identify 
the corresponding $j$-coordinates of its peaks 
with parts of partitions described in the Lemma (see Figure~\ref{fig:AXEL}).

\noindent
Hence, we can reformulate this Lemma as
\begin{equation}
q^{t^2}\qBin{L-t}{t}{q}=
\sum_{\mathbf{P}_{0,0}^t(0,L)}q^{\sum_{j=1}^{L-1}j\sigma_j},
\label{eq:2.1}
\end{equation}
where the symbol $\Sigma_{\mathbf{P}_{s,b}^t(i,f)}$ denotes 
the sum over all admissible paths
$\in\mathbf{P}_{s,b}^t(i,f)$. The advantage of the path interpretation is that 
while partitions always have non-negative parts (by definition), $j$-coordinates of peaks 
in general may take on negative values.

If we move $\mathbf{P}_{0,0}^t(0,L)$ to the left by $M$ units, then the resulting 
path space $\mathbf{P}_{0,0}^t(-M,L-M)$ can be used to prove that    
\begin{equation}
q^{t^2-Mt}\qBin{L-t}{t}{q}=
\sum_{\mathbf{P}_{0,0}^t(-M,L-M)}q^{\sum_{j=1-M}^{L-M-1}j\sigma_j}.
\label{eq:2.2}
\end{equation}
Indeed, using \eqn{2.1} we have
\begin{align}
\sum_{\mathbf{P}_{0,0}^t(-M,L-M)}q^{\sum_{j=1-M}^{L-M-1}j\sigma_j}&=
\sum_{\mathbf{P}_{0,0}^t(0,L)}q^{\sum_{j=1}^{L-1}(j-M)\sigma_j}=
q^{-Mt}\sum_{\mathbf{P}_{0,0}^t(0,L)}q^{\sum_{j=1}^{L-1}j\sigma_j}\notag\\
&=q^{t^2-Mt}\qBin{L-t}{t}{q}.
\label{eq:2.3}
\end{align}
And so,
\begin{equation}
\sum_{t\ge0}q^{t^2-Mt}\qBin{L-t}{t}{q}=
\sum_{\mathbf{P}_{0,0}(-M,L-M)}q^{\sum_{j=1-M}^{L-M-1}j\sigma_j},
\label{eq:2.4}
\end{equation}
where $\mathbf{P}_{s,b}(i,f)$ is defined the same way as 
$\mathbf{P}_{s,b}^t(i,f)$, except that we no longer require that the 
number of peaks is exactly $t$.

More generally, one can easily show that for $s,b\in\{0,1\}$
\begin{equation}
f_{s,b}(L,M,q)=C_{s,b}(-M,L-M,q)
\label{eq:2.5}
\end{equation}
with
\begin{equation}
f_{s,b}(L,M,q):=\sum_{t\ge0}q^{t^2+st-Mt}
\qBin{L-t-s-b}{t}{q}
\label{eq:2.6}
\end{equation}
and
\begin{equation}
C_{s,b}(i,f,q):=\sum_{\mathbf{P}_{s,b}(i,f)}q^{\sum_{j=i+1}^{f-1}j\sigma_j}.
\label{eq:2.7}
\end{equation}
Next, for $0\le M\le L$ every admissible path $\in\mathbf{P}_{s,b}(-M,L-M)$
consists of two pieces joined together at point $(0;s'=0,1)$. The first
piece belongs to $\mathbf{P}_{s,s'}(-M,0)$ and the second one to 
$\mathbf{P}_{s',b}(0,L-M)$. This observation is equivalent to
\begin{equation}
C_{s,b}(-M,L-M,q)=\sum_{s'=0}^1 C_{s,s'}(-M,0,q)C_{s',b}(0,L-M,q).
\label{eq:2.8}
\end{equation}
Now, because
\begin{equation}
\sum_{\mathbf{P}_{s,s'}(-M,0)}q^{\sum_{j=1-M}^{-1}j\sigma_j}=
\sum_{\mathbf{P}_{s',s}(0,M)}(\frac{1}{q})^{\sum_{j=1}^{M-1}j\sigma_j}
\label{eq:2.9}
\end{equation}
we infer that
\begin{equation}
C_{s,s'}(-M,0,q)=C_{s',s}(0,M,\frac{1}{q}).
\label{eq:2.10}
\end{equation}
Next, combining \eqn{2.5}, \eqn{2.8} and \eqn{2.10}, we arrive at
\begin{equation}
f_{s,b}(L,M,q)=\sum_{s'=0}^1 f_{s',s}(M,0,\frac{1}{q})f_{s',b}(L-M,0,q).
\label{eq:2.11}
\end{equation}
The desired formula \eqn{1.14} is an easy consequence of \eqn{2.11} with $s=b=0$ 
and the Rogers-Ramanujan identities \eqn{1.3} and \eqn{1.4}, which we restate again as
\begin{equation}
f_{s,0}(L,0,q)=\sum_{t\ge0}q^{t^2+s t}\qBin{L-t-s}{t}{q}=
\begin{cases}e_L(q), & \mbox{if } s=0, \\
d_L(q), & \mbox{if } s=1. \end{cases}
\label{eq:2.12}
\end{equation}
Remark: If we set $s=1,b=0,q=1$ in \eqn{2.11}, we immediately derive the following
well known identity for the Fibonacci numbers $\mathrm{Fi}(L)$:
\begin{equation}
\mathrm{Fi}(L)=\mathrm{Fi}(M)\mathrm{Fi}(L-M+1)+\mathrm{Fi}(M-1)\mathrm{Fi}(L-M).
\label{eq:2.13}
\end{equation}
If we perform the substitution $M\rightarrow-M$ and use,
according to \eqn{1.10},
\begin{equation}
\mathrm{Fi}(-M)=(-1)^{M+1}\mathrm{Fi}(M)
\label{eq:2.14}
\end{equation}
in \eqn{2.13}, we obtain
\begin{equation}
(-1)^M \mathrm{Fi}(L)=\mathrm{Fi}(M+1)\mathrm{Fi}(L+M)-\mathrm{Fi}(M)\mathrm{Fi}(L+M+1),
\label{eq:2.15}
\end{equation}
which is a specialization of the Euler-Cassini formula. We would like to point out that
\eqn{2.13} is ``minus sign" free. As a result, the combinatorial proof of \eqn{2.13} given 
here is very different from that of Werman and Zeilberger \cite{WZ}. 
Namely, their proof of \eqn{2.15}
given in \cite{WZ} made essential use of involution technique.

We conclude this section by pointing out that our analysis can be trivially 
extended to show that
\begin{align}
f_{s,b}(L,M,q)&=\sum_{s'=0}^1 f_{s,s'}(M+x,M,q)\;q^{s'x}f_{s',b}(L-M-x,-x,q)\notag\\
&=\sum_{s'=0}^1 q^{s'x}f_{s',s}(M+x,x,\frac{1}{q})\;f_{s',b}(L-M-x,-x).
\label{eq:2.16}
\end{align}
Note that \eqn{2.11} is \eqn{2.16} with $x=0$.

\section{Variants of the Andrews-Gordon Identities}\label{sec:3}

For $\nu\in\mathbb{Z}_{>0}$, the analytical generalizations of the Rogers-Ramanujan 
identities known as Andrews-Gordon identities \cite{A2} can be stated as
\begin{align}
\sum_{n_1,n_2,\ldots,n_\nu}&\frac{q^{(N_1^2+N_2^2+\cdots+N_\nu^2)+(N_s+N_{s+1}+\cdots+N_\nu)}}
{(q)_{n_1}(q)_{n_2}\ldots(q)_{n_\nu}} \notag\\
& = \frac{1}{(q)_\infty}\sum_{j=-\infty}^\infty 
\left\{q^{2(2\nu+3)j^2+j(2\nu+3-2s)}-q^{(2j+1)((2\nu+3)j+s)}\right\}\notag\\
& = \frac{1}{\prod_{n\not\equiv 0,\pm s\,(\mymod 2\nu+3)}(1-q^n)}, \quad s=1,2,\ldots,\nu+1
\label{eq:3.1}
\end{align}
with
\begin{equation}
N_i=\begin{cases}n_i+n_{i+1}+\cdots+n_\nu, & \mbox{if } 1\le i\le\nu, \\
0, & \mbox{if } i=\nu+1.\end{cases}
\label{eq:3.2}
\end{equation}
Here and throughout, we adopt the convention that in the product 
\[
  \prod_{n\not\equiv 0,\pm s \,(\mymod 2\nu+3)}(1-q^n)
\]
$n$ takes on positive integer values not
congruent to $0,\pm s \mbox{ modulo }2\nu+3$. Clearly, if $\nu=1$, \eqn{3.1} 
reduces to \eqn{1.1}.
As in the case of the Rogers-Ramanujan identities, 
the identities \eqn{3.1} have polynomial analogs.
To describe these polynomial versions we need to introduce polynomials $\tilde{F}_{s,b}(L,q)$
defined for $0\le s,b\le\nu$ as follows
\begin{equation}
\tilde{F}_{s,b}(L,q):=\sum_{\mathbf{n}}q^{(N_1^2+\cdots+N_\nu^2)+(N_{s+1}+\cdots+N_\nu)}
\qBin{\mathbf{n}+\mathbf{m}}{\mathbf{n}}{q}
\label{eq:3.3}
\end{equation}
with
\begin{equation}
\mathbf{n}=(n_1,n_2,\ldots,n_\nu),\;\;\;\mathbf{m}=(m_1,m_2,\ldots,m_\nu)
\label{eq:3.4}
\end{equation}
and
\begin{equation}
\qBin{\mathbf{n}+\mathbf{m}}{\mathbf{n}}{q}=\prod_{i=1}^\nu
\qBin{n_i+m_i}{n_i}{q},
\label{eq:3.5}
\end{equation}
where
\begin{equation}
m_i=L-2(N_1+N_2+\cdots+N_i)-\chi(i>s)(i-s)-\chi(i>b)(i-b)
\label{eq:3.6}
\end{equation}
and
\begin{equation}
\chi(i>a)=\begin{cases}1, & \mbox{if } i>a, \\
0, & \mbox{if } i\le a.\end{cases}
\label{eq:3.7}
\end{equation}
Next, for $L\equiv s+b \,(\mymod 2)$ and $1\le s,b\le\nu+1$, we define polynomials $B_{s,b}(L,q)$ as
\begin{align}
B_{s,b}(L,q):=\sum_{j=-\infty}^\infty
\Bigg\{ & q^{2(2\nu+3)j^2+j(2\nu+3-2s)}
\qBin{L}{\frac{L+s-b}{2}-j(2\nu+3)}{q} \notag\\
 & -q^{(2j+1)((2\nu+3)j+s)}
\qBin{L}{\frac{L-s-b}{2}-j(2\nu+3)}{q} \Bigg\}.
\label{eq:3.8}
\end{align}
Equipped with these definitions we are in the position to state the 
polynomial analogs of \eqn{3.1}, namely
\begin{equation}
\tilde{F}_{s,b}(L,q)=\begin{cases}
B_{s+1,b+1}(L,q), & \mbox{if } L \equiv s+b \, (\mymod 2), \\
B_{(2\nu+3)-(s+1),b+1}(L,q), & \mbox{if } L \not\equiv s+b \, (\mymod 2), \end{cases}
\label{eq:3.9}
\end{equation}
with $0\le s,b\le\nu$.

For $b=\nu$, formulas \eqn{3.9} first appeared in the works of 
Foda, Quano \cite{FQ} and Kirillov \cite{Kir}, for other values of $b$, these formulas were 
derived in \cite{BMS}. It is important to keep in mind 
that in case $s\neq\nu$ and $b\neq\nu$, 
the summands in \eqn{3.3} may be non-zero in value even if $n_\nu=-1$.
 
To prove \eqn{3.9} the authors of \cite{BMS} showed that both sides of \eqn{3.9} 
satisfy identical recurrences for $1\le b\le\nu$,
\begin{equation}
\begin{cases}\tilde{F}_{s,0}(L,q)=\tilde{F}_{s,1}(L-1,q),\\
\tilde{F}_{s,b}(L,q)=\tilde{F}_{s,b-1}(L-1,q)+\tilde{F}_{s,b+1-\delta_{b,\nu}}(L-1,q)
+(q^{L-1}-1)\tilde{F}_{s,b}(L-2,q)\end{cases}
\label{eq:3.10}
\end{equation}
and the initial conditions
\begin{equation}
\tilde{F}_{s,b}(0,q)=\delta_{s,b}
\label{eq:3.11}
\end{equation}
where the Kronecker delta function $\delta_{i,j}$ is defined as usual as
\begin{equation}
\delta_{i,j}=\begin{cases}1, & \mbox{if } i=j, \\
0, & \mbox{otherwise.} \end{cases}
\label{eq:3.12}
\end{equation}
We would like to emphasize that there is more than one way to finitize the Andrews-Gordon
identities. In particular, Warnaar \cite{W} found polynomial versions of \eqn{3.1} involving 
$q$-multinomial coefficients.

Now we would like to alter the 
Andrews-Gordon identities in the spirit of Garrett et al.~\cite{GIS}.
However, because we have more than one summation variable, it is not immediately 
clear how to accomplish this. 
Our guiding principle is that additional linear terms should modify 
the recurrences \eqn{3.10}
in a minimal way by simple shifts, since this is precisely what happened in case of 
Rogers-Ramanujan identities. More specifically, we would like to have polynomials
$F_{s,b}(L,M,q)$ satisfying the following relations for $1 \le b\le\nu,\;\;0\le s\le\nu$
\begin{align}
\left\{\begin{array}{l}F_{s,0}(L,M,q)=F_{s,1}(L-1,M,q) \\
F_{s,b}(L,M,q)=F_{s,b-1}(L-1,M,q)+F_{s,b+1-\delta_{b,\nu}}(L-1,M,q)\\
\qquad\qquad\quad\:\:\: {} +(q^{L-M-1}-1)F_{s,b}(L-2,M,q).\end{array}\right.
\label{eq:3.13}
\end{align}
The above requirement leads us to define
\begin{equation}
F_{s,b}(L,M,q):=\sum_{\mathbf{n}}q^{(N_1^2+N_2^2+\cdots+N_\nu^2)+(N_{s+1}+\cdots+N_\nu)-MN_1}
\qBin{\mathbf{n}+\mathbf{m}}{\mathbf{n}}{q}
\label{eq:3.14}
\end{equation}
with $s,b\in\{0,1,\ldots,\nu\}$ and the rest of notations the same as in \eqn{3.3}.
In the Appendix we will prove that these polynomials
indeed satisfy the recurrences \eqn{3.13}.

Now, since $F_{s,b}(L,M,q)$ and $F_{s,b}(L-M,0,q)=\tilde{F}_{s,b}(L-M,q)$ satisfy the same
recursion relations we can write
\begin{equation}
F_{s,b}(L,M,q)=\sum_{s'=0}^\nu A_{s,s'}(M,q)\tilde{F}_{s',b}(L-M,q).
\label{eq:3.15}
\end{equation}
The connection coefficients $A_{s,s'}(M,q)$ can be easily determined from the boundary 
conditions
\begin{equation}
F_{s,b}(L,L,q)=\sum_{s'=0}^\nu A_{s,s'}(L,q)\tilde{F}_{s',b}(0,q)
\stackrel{\mbox{\scriptsize{by \eqn{3.11}}}}{=}
\sum_{s'=0}^\nu A_{s,s'}(L,q)\delta_{s',b}=A_{s,b}(L,q).
\label{eq:3.16}
\end{equation}
Making use of
\begin{equation}
\qBin{n+m}{n}{\frac{1}{q}}=q^{-nm}
\qBin{n+m}{n}{q}
\label{eq:3.17}
\end{equation}
one can easily verify that
\begin{equation}
F_{s,b}(L,L,q)=F_{b,s}(L,0,\frac{1}{q})=\tilde{F}_{b,s}(L,\frac{1}{q}).
\label{eq:3.18}
\end{equation}
Hence,
\begin{equation}
F_{s,b}(L,M,q)=\sum_{s'=0}^\nu\tilde{F}_{s',s}(M,\frac{1}{q})\tilde{F}_{s',b}(L-M,q),
\label{eq:3.19}
\end{equation}
which is a perfect analog of formula \eqn{2.11}. Recalling \eqn{3.9}, we can rewrite 
\eqn{3.19} as
\begin{align}
F_{s,b}(L,M,q)&=\sum_{\stackrel{s'=0}{s+s'\equiv M\,(\mymod 2)}}^\nu B_{s'+1,s+1}(M,\frac{1}{q})
\tilde{F}_{s',b}(L-M,q)\notag\\
&+\sum_{\stackrel{s'=0}{s+s'\not\equiv M\,(\mymod 2)}}^\nu 
B_{(2\nu+3)-(s'+1),s+1}(M,\frac{1}{q})\tilde{F}_{s',b}(L-M,q).
\label{eq:3.20}
\end{align}
In the limit $L\rightarrow\infty$, \eqn{3.20} gives
\begin{align}
\sum_{\mathbf{n}} &\frac{q^{N_1^2+N_2^2+\cdots+N_\nu^2+(N_s+\cdots+N_\nu)-MN_1}}
{(q)_{n_1}(q)_{n_2}\ldots(q)_{n_\nu}}\notag\\
&= \sum_{\stackrel{s'=1}{s+s'\equiv M\,(\mymod 2)}}^{\nu+1}\frac{B_{s',s}(M,\frac{1}{q})}
{\prod_{n\not\equiv 0,\pm s'\,(\mymod 2\nu+3)}(1-q^n)}\notag\\
&+\sum_{\stackrel{s'=1}{s+s'\not\equiv M\,(\mymod 2)}}^{\nu+1}\frac{B_{2\nu+3-s',s}(M,\frac{1}{q})}
{\prod_{n\not\equiv 0,\pm s'\,(\mymod 2\nu+3)}(1-q^n)},
\label{eq:3.21}
\end{align}
where we used the limiting formulas
\begin{equation}
\lim_{L\rightarrow\infty}\tilde{F}_{s-1,b}(L,q)=\frac{1}
{\prod_{n\not\equiv 0,\pm s\,(\mymod 2\nu+3)}(1-q^n)},
\label{eq:3.22}
\end{equation}
which follow from 
\begin{equation}
\lim_{L\rightarrow\infty}\qBin{L}{n}{q}=\frac{1}{(q)_n}
\label{eq:3.23}
\end{equation}
and \eqn{3.1}. It is easy to check that in case $\nu=1,s=2$ \eqn{3.21} reduces to \eqn{1.12}.

\section{Concluding Remarks}\label{sec:4}

The interested reader may wonder if the combinatorial analysis given in Section~\ref{sec:2}
can be upgraded to explain the formulas \eqn{3.19}. The answer to this
question is affirmative. 
However,
for $\nu>1$ the path interpretation of the $F_{s,b}(L,M,q)$
polynomials is much 
more involved than that of the $f_{s,b}(L,M,q)$ polynomials considered 
in Section~\ref{sec:2}. Here, one should deal with 
peaks of different heights \cite{Bre} and in addition with certain
boundary 
defects. We plan to 
come back to the combinatorial derivation of \eqn{3.19} in our future work. Here, we confine 
ourselves to  
remark that the introduction of an additional linear term $-MN_1$
in \eqn{3.21} amounts to the shift to the left by $M$ units of Bressoud's path 
described in \cite{Bre}.

\medskip
However, with respect to peaks with different heights we want to
mention that these pop-up also in connection with another 
polynomial version of an identity of Garrett-Ismail-Stanton type.
Namely, for integers $L,m \geq 0$ one has
\begin{equation}
\label{P1}
q^m (q;q^2)_m \sum_{t\geq 0} \qBin{L}{2t+1}{q} q^{2 t^2 +
  2(m+1)t} =
S_m(q) T_{L+m}(q) - T_m(q) S_{L+m}(q),
\end{equation}
where 
\begin{equation*}
S_m(q) = \sum_{t \geq 0} q^{2 t^2} \qBin{L}{2t}{q}
\end{equation*}
and
\begin{equation*}
T_m(q) = \sum_{t \geq 0}q^{2 t^2+2 t} \qBin{L}{2t+1}{q}
\end{equation*}
are the Andrews-Santos polynomials discussed in \cite{AS}.
Identity \eqref{P1} arose in work of Andrews et al.~\cite{AKPZ}
and yields, in the limit $L \rightarrow \infty$, a combination
of Slater's identities (38) and (39) from \cite{Slater}. 

The methods of Sections~\ref{sec:1} and \ref{sec:2} above can be applied and lead to
the following generalization which extends \eqref{P1} also to
negative integers: For $L\geq 0$ and arbitrary integer $M$,
\begin{equation}
\label{P2}
\sum_{t\geq 0} \qBin{L}{2t+1}{q} q^{2 t^2 - 2 M t} =
q^{M+1} S_{M+1}(\frac{1}{q}) T_{L-M-1}(q) +
q^M T_{M+1}(\frac{1}{q}) S_{L-M-1}(q).
\end{equation}
Here we understand that for negative indices, i.e., for $m\geq 0$,
one has
\begin{equation*}
S_{-m}(q) = (-1)^m \frac{q^{m^2}}{(q;q^2)_m} S_m(\frac{1}{q})
\end{equation*}
and 
\begin{equation*}
T_{-m}(q) = (-1)^{m+1} \frac{q^{m^2-1}}{(q;q^2)_m} T_m(\frac{1}{q}).
\end{equation*}
In the limit $L\rightarrow \infty$ one obtains from \eqref{P2} 
another new identity of Garrett-Ismail-Stanton type. The proof
and the underlying combinatorics will be presented in a 
forthcoming paper.

\medskip
Our variants of the Andrews-Gordon identities were determined by 
the polynomial versions \eqn{3.9}
and a requirement that the introduction of additional 
linear terms should modify the recursion
relations \eqn{3.10} by trivial shifts as in \eqn{3.13}. We intend to use more general
polynomial versions of Andrews-Gordon identities containing $\nu$ finitization parameters
to investigate the most general multisum
\begin{equation}
\sum_{\mathbf{n}}\frac{q^{N_1^2+\cdots+N_\nu^2-M_1N_1-M_2N_2-MN_\nu}}
{(q)_{n_1}(q)_{n_2}\ldots(q)_{n_\nu}}.
\label{eq:4.1}
\end{equation}

Finally, we would like to mention that many new generalizations of Rogers-Ramanujan 
identities were introduced in \cite{BM} and proven in \cite{BMS} and \cite{FLW}. 
Techniques developed in this paper are adequate to produce and prove variants of all 
these identities.

\subsection*{Acknowledgment}

We are grateful to Prof.~G.E.~Andrews, Prof.~J.~Cigler and Prof.~C.~Krattenthaler 
for their interest and stimulating discussions. Special thanks to Dr.~A.~Riese
for his help in preparing this manuscript. This paper was completed while 
the first author was visiting RISC--Linz, Austria. He would like to thank faculty and 
staff for their warm hospitality.


\section{Appendix}

Here we will prove the recurrences \eqn{3.13}. To this end we need to define vectors
${\mathbf{e}}_i$ and ${\mathbf{E}}_{a,b}$ as
\begin{equation}
(\mathbf{e}_i)_j=\begin{cases}1, & \mbox{if }i=j \mbox{ and } 1\le i\le\nu, \\
0, & \mbox{otherwise,}\end{cases}
\label{eq:A.1}
\end{equation}
and
\begin{equation}
\mathbf{E}_{a,b}=\sum_{i=a}^b\mathbf{e}_i.
\label{eq:A.2}
\end{equation}
The first recurrence in \eqn{3.13} is trivial. To prove the second relation in 
\eqn{3.13}, we expand $F_{s,b}(L,M,q)$ in a telescopic fashion as
\begin{align}
F_{s,b}(L,M,q)&=\sum_{\mathbf{n}}q^{\Phi_s(\mathbf{N},M)}
\qBin{\mathbf{n}+\mathbf{m}-\mathbf{E}_{1,b}}{\mathbf{n}}{q} \notag\\
&+\sum_{\mathbf{n}}q^{\Phi_s(\mathbf{N},M)}
\qBin{\mathbf{n}+\mathbf{m}-\mathbf{E}_{1,b}}{\mathbf{n}-\mathbf{e}_b}{q} q^{m_b} \notag\\
&+\sum_{\mathbf{n}}q^{\Phi_s(\mathbf{N},M)}
\qBin{\mathbf{n}+\mathbf{m}-\mathbf{E}_{1,b-1}}{\mathbf{n}-\mathbf{e}_{b-1}}{q} q^{m_{b-1}}\notag\\
&+\sum_{\mathbf{n}}q^{\Phi_s(\mathbf{N},M)}
\qBin{\mathbf{n}+\mathbf{m}-\mathbf{E}_{1,b-2}}{\mathbf{n}-\mathbf{e}_{b-2}}{q} q^{m_{b-2}}\notag\\
& \ldots \notag \\
&+\sum_{\mathbf{n}}q^{\phi_s(\mathbf{N},M)}
\qBin{\mathbf{n}+\mathbf{m}-\mathbf{E}_{1,1}}{\mathbf{n}-\mathbf{e}_1}{q} q^{m_1}, 
\label{eq:A.3}
\end{align}
where $\mathbf{N}=(N_1,N_2,\ldots,N_\nu)$ and
\begin{equation}
\Phi_s(\mathbf{N},M)=(N_1^2+N_1^2+\cdots+N_\nu^2)+(N_{s+1}+\cdots+N_\nu)-MN_1,
\label{eq:A.4}
\end{equation}
with the rest of notations the same as in \eqn{3.3}.

It is important to remember that the vectors $\mathbf{m}$ in \eqn{A.3} are actually functions of 
$\mathbf{N}$ and $L$ as can be seen from \eqn{3.6}.
To make sure that \eqn{A.3} is a correct expansion of $F_{s,b}(L,M,q)$ we merge the first 
and second sum in \eqn{A.3} into a single sum using
\begin{equation}
\qBin{n+m}{n}{q}=
\qBin{n+m-1}{n}{q}+q^m
\qBin{n-1+m}{n-1}{q}.
\label{eq:A.5}
\end{equation}
This single sum, in turn, can be merged with the third sum in \eqn{A.3}. This process can be
repeated until all sums are merged together to yield $F_{s,b}(L,M,q)$. It is trivial to
recognize the first sum in \eqn{A.3} as $F_{s,b+1-\delta_{b,\nu}}(L-1,M,q)$. With regard
to the last sum in \eqn{A.3}, we perform the change $n_1\rightarrow n_1+1$ to recognize
it as $q^{L-M-1}F_{s,b}(L-2,M,q)$. So, it follows that
\begin{align}
&F_{s,b}(L,M,q)-F_{s,b+1-\delta_{b,\nu}}(L-1,M,q)-q^{L-M-1}F_{s,b}(L-2,M,q) \notag\\
&\qquad=\sum_{i=2}^b\sum_{\mathbf{n}}q^{\Phi_s(\mathbf{N},M)+m_i}
\qBin{\mathbf{n}+\mathbf{m}-\mathbf{E}_{1,i}}{\mathbf{n}-\mathbf{e}_i}{q}.
\label{eq:A.6}
\end{align}
If $b=1$, then the rhs of \eqn{A.6} is just zero, so
\begin{equation}
F_{s,1}(L,M,q)=F_{s,2}(L-1,M,q)+q^{L-M-1}F_{s,1}(L-2,M,q).
\label{eq:A.7}
\end{equation}
Combining \eqn{A.7} with the first recurrence in \eqn{3.13}
\begin{equation}
F_{s,0}(L-1,M,q)=F_{s,1}(L-2,M,q)
\label{eq:A.8}
\end{equation}
we obtain
\begin{equation}
F_{s,1}(L,M,q)=F_{s,0}(L-1,M,q)+F_{s,2}(L-1,M,q)+(q^{L-M-1}-1)F_{s,1}(L-2,M,q),
\label{eq:A.9}
\end{equation}
as desired.

If $b\neq 0,1$, we perform an $``i"$ dependent change of the summation variables in~\eqn{A.6}
\begin{equation}
\mathbf{n}\rightarrow\mathbf{n}+
(\mathbf{e}_i-\mathbf{e}_{i-1})-(\mathbf{e}_b-\mathbf{e}_{b-1})
\label{eq:A.10}
\end{equation}
to obtain
\begin{align}
& F_{s,b}(L,M,q)-F_{s,b+1-\delta_{b,\nu}}(L-1,M,q)-q^{L-M-1}F_{s,b}(L-2,M,q) \notag\\
& \qquad = \sum_{i=1}^{b-1}\sum_{\mathbf{n}}q^{\Phi_s(\mathbf{N},M)+m_i+m_b-m_{b-1}}
\qBin{\mathbf{n}+\mathbf{m}-\mathbf{E}_{1,b-1}+
\mathbf{e}_{b-1}-\mathbf{e}_b-\mathbf{E}_{i,b-1}}
{\mathbf{n}+\mathbf{e}_{b-1}-\mathbf{e}_b-\mathbf{e}_i}{q}.
\label{eq:A.11}
\end{align}
Next, we rewrite the polynomial $F_{s,b-1}(L-1,M,q)$ in terms of the same $\mathbf{n}$
vectors as in \eqn{A.11}
\begin{equation}
F_{s,b-1}(L-1,M,q)=\sum_{\mathbf{n}}q^{\Phi_s(\mathbf{N},M)+m_b-m_{b-1}+1}
\qBin{\mathbf{n}+\mathbf{m}-\mathbf{E}_{1,b-1}+
\mathbf{e}_{b-1}-\mathbf{e}_b}{\mathbf{n}+\mathbf{e}_{b-1}-\mathbf{e}_b}{q}
\label{eq:A.12}
\end{equation}
and then expand it in a telescopic fashion to get
\begin{align}
&F_{s,b-1}(L-1,M,q)\notag\\
&\qquad= \sum_{\mathbf{n}}q^{\Phi_s(\mathbf{N},M)+m_b-m_{b-1}+1}
\qBin{\mathbf{n}+\mathbf{m}-\mathbf{E}_{1,b-1}+
\mathbf{e}_{b-1}-\mathbf{e}_b-\mathbf{E}_{1,b-1}}
{\mathbf{n}+\mathbf{e}_{b-1}-\mathbf{e}_b}{q} \notag\\
& \qquad + \sum_{i=1}^{b-1}\sum_{\mathbf{n}}q^{\Phi_s(\mathbf{N},M)+m_i+m_b-m_{b-1}}
\qBin{\mathbf{n}+\mathbf{m}-\mathbf{E}_{1,b-1}+
\mathbf{e}_{b-1}-\mathbf{e}_b-\mathbf{E}_{i,b-1}}
{\mathbf{n}+\mathbf{e}_{b-1}-\mathbf{e}_b-\mathbf{e}_i}{q}.
\label{eq:A.13}
\end{align}
Once again, performing the change 
$\mathbf{n}\rightarrow\mathbf{n}-\mathbf{e}_{b-1}+\mathbf{e}_b$ in the first sum
on the rhs of \eqn{A.13} we recognize it as $F_{s,b}(L-2,M,q)$.  
Hence,
\begin{align}
&F_{s,b-1}(L-1,M,q)-F_{s,b}(L-2,M,q)\notag\\
& \qquad =\sum_{i=1}^{b-1}\sum_{\mathbf{n}}q^{\Phi_s(\mathbf{N},M)+m_i+m_b-m_{b-1}}
\qBin{\mathbf{n}+\mathbf{m}+
\mathbf{e}_{b-1}-\mathbf{e}_b-\mathbf{E}_{i,b-1}}
{\mathbf{n}+\mathbf{e}_{b-1}-\mathbf{e}_b-\mathbf{e}_i}{q}.
\label{eq:A.14}
\end{align}
Comparing the rhs of \eqn{A.11} and \eqn{A.14} we immediately infer that
\begin{align}
F_{s,b}(L,M,q) & = F_{s,b-1}(L-1,M,q)+F_{s,b+1-\delta_{b,\nu}}(L-1,M,q)\notag\\
& + (q^{L-M-1}-1)F_{s,b}(L-2,M,q),
\label{eq:A.15}
\end{align}
as desired.

\end{document}